\documentclass[12pt]{article}
\def\today{January 14, 2003}
\usepackage{amsmath}
\usepackage{amssymb}
\usepackage{amsthm}
\usepackage[pdfpagemode=None]{hyperref}
\usepackage{psfig}
\usepackage{epsf}
\usepackage{psfrag}
\newtheorem{thm}{Theorem}[section]
\newtheorem{co}[thm]{Corollary}

\newtheorem{assumption}[thm]{Assumption}

\newtheorem{definition}[thm]{Definition}

\newtheorem{example}[thm]{Example}

\newtheorem{remark}[thm]{Remark}

\newtheorem{tab}{Table}

\newcommand{\Section}[1]{\section{#1}\setcounter{equation}{0}}

\newcommand{\eqr}[1]{~\mbox{$(${\rm \ref{#1}}$)$}}

\newcommand{\diag}{{\rm diag}\,}

\newcommand{\V}{\mathcal{V}}

\title{Generalized PSK in Space Time Coding\footnote{The author is supported
    in part by NSF grants DMS-00-72383 and CCR-02-05310 and by a fellowship from the
    Center of Applied Mathematics at the University of Notre
    Dame.}}
\date{{\normalsize \today}}%

\author{Guangyue Han\\
  {\normalsize Department of Mathematics}\vspace{-1mm} \\
  {\normalsize University of Notre Dame}\vspace{-1mm} \\
  {\normalsize Notre Dame, IN 46556.}\\
  {\normalsize {\em e-mail:\/} Han.13@nd.edu}\vspace{-1mm}\\
  {\normalsize {\em URL:} http://www.nd.edu/\~{}eecoding/} }

\begin{document}\maketitle\thispagestyle{empty}

\begin{abstract}
  A wireless communication system using multiple antennas
  promises reliable transmission under Rayleigh flat fading
  assumptions.  Design criteria and practical schemes have been
  presented for both coherent and non-coherent communication
  channels. In this paper we generalize one dimensional phase
  shift keying (PSK) signals and introduce space time constellations from generalized
  phase shift keying (GPSK) signals based on the complex and real orthogonal designs. The
  resulting space time constellations reallocate the energy for
  each transmitting antenna and feature good diversity products,
  consequently their performances are better than some of the
  existing comparable codes. Moreover since the maximum likelihood (ML) decoding of our
  proposed codes can be decomposed to one dimensional PSK signal
  demodulation, the ML decoding of our codes can be
  implemented in a very efficient way.
\end{abstract}


\Section{Introduction and Model}

Multiple antenna communication systems have been actively studied
recently~\cite{al98,ho00,ho00a,ta98}. By exploiting the temporal
and spatial diversity at both transmitter side and receiver side,
such kind of systems can increase the channel capacity compared to
single antenna communication systems, consequently promising more
reliable data transmission for high rate applications.

We investigate a communication system under Rayleigh flat fading
assumptions and we further assume fading is quasi-static over a
time period of length $T$. Denote by $M$ or $N$ the number of
transmitting or receiving antennas, respectively. Let $\rho$
represent the expected signal-to-noise ratio (SNR) at each
receiving antenna. For the system above, the basic equation
between the received signal $R$, which is a $T \times N$ matrix,
and the transmitted signal $\Phi$, which is chosen from a $T
\times M$ matrix constellation $\V=\{ \Phi_1, \Phi_2, \cdots,
\Phi_L\}$ ($L$ is the constellation size), is given through:
$$
R=\sqrt{\frac{\rho}{M}}\Phi H+W,
$$
where the $M \times N$ matrix $H$ accounts for the
multiplicative complex Gaussian fading coefficients and the $T
\times N$ matrix $W$ accounts for the additive white Gaussian
noise. The entries $h_{m,n}$ of the matrix $H$ as well as the
entries $w_{t,n}$ of the matrix $W$ are assumed to have a
statistically independent complex normal distribution
$\mathcal{CN}(0,1)$. One can verify that the transmission rate is
determined by $L$ and $T$:
$$
\mathtt{R}=\frac{\log_2(L)}{T}.
$$

When the fading coefficients are unknown to the transmitter
however known to the receiver, it is proven~\cite{fo96a1,te99}
that the channel capacity will increase linearly with $\min \{M,
N\}$. Since $H$ is known to the receiver, for a received signal
$R$ the maximum likelihood (ML) decoder will take the following
evaluation to resolve the most likely sent signal (codeword):
$$
\hat{\Phi}=\arg \min_{\Phi \in \V}
{\left\|R-\sqrt{\frac{\rho}{M}} \Phi H \right\|},
$$
where $\| \cdot \|$ represents Frobenius norm. Let $P_{\Phi_l,
\Phi_{l'}}$ denote the probability that ML decoder is mistaking
$\Phi_l$ for $\Phi_{l'}$. The upper bound of this error
probability has been derived in~\cite{ta98}:
$$
P_{\Phi_l, \Phi_{l'}} \leq \frac{1}{2} \prod_{m=1}^M \left[
1+\frac{\rho T}{4M} \delta_m^2(\Phi_l-\Phi_{l'})\right]^{-N}.
$$

Through the analysis of the above upper bound, design criteria
for designing codes for the coherent channel with ideal channel
state information (CSI) are proposed in~\cite{ta98}:

\textit{The Rank Criterion}: In order to achieve the maximum
diversity $MN$, the matrix $B(\Phi_l, \Phi_l')=\Phi_l-\Phi_{l'}$
has to have full rank for any codewords $\Phi_l$ and $\Phi_{l'}$.
If $B(\Phi_l, \Phi_{l'})$ has minimum rank $r$ over the set of two
tuples of distinct codewords, then a diversity $rM$ is achieved.

\textit{The Determinant Criterion}: Suppose that a diversity
benefit of $rM$ is our target. The minimum of $r$th roots of the
sum of determinants of all $r \times r$ principal cofactors of
$A(\Phi_l, \Phi_{l'})=B(\Phi_l, \Phi_{l'}) B^*(\Phi_l, \Phi_{l'})$
taken over all pairs of distinct codewords $\Phi_l$ and
$\Phi_{l'}$ corresponds to the coding advantage, where $r$ is the
rank of $A(\Phi_l, \Phi_{l'})$.

Note that for the case when $T=M$ and full diversity is achieved,
the above design criteria can be simplified as follows:

Construct a constellation of matrix $\V=\{\Phi_1, \Phi_2, \cdots,
\Phi_L\}$ such that the diversity product~\cite{ho00}
$$
\prod \V=\min_{l \neq l'} \frac{1}{2}
{|\det(\Phi_l-\Phi_{l'})|}^{1/M}
$$
is as large as possible.

Orthogonal designs have been investigated for the constellation
construction for coherent channels. Using a complex orthogonal
design, Alamouti~\cite{al98} proposed a very simple transmitter
diversity scheme with $2$ transmitting antennas. A very
interesting coding scheme from the real orthogonal designs is
presented in~\cite{ta99}. Applying the similar idea as
in~\cite{al98}, the authors also explain how the symmetric
structure of the orthogonal design codes leads to a much simpler
decoding algorithm. However no explicit constellations have been
studied in details in these work.

For a slow fading channel, for instance a fixed wireless
communication system, the coherent assumption is reasonable
because the transmitter can send pilot signals which allow the
receiver to estimate the fading coefficients accurately. However
in certain situations, due to the difficulty of measuring the
fading coefficients under the practical communication environment
(for instance the limited resources or a fast fading link) the
above coherent model would be questionable. Because of this,
researchers put their efforts into the case when the channel is
non-coherent, i.e., the CSI is not known by either the transmitter
or the receiver.

In~\cite{ho00a}, Hochwald and Marzetta study unitary space-time
modulation for a non-coherent channel. We will use the same
notations as in the coherent case, so the basic equation will be
still the same:
$$
R=\sqrt{\frac{\rho}{M}}\Phi H+W,
$$
however in non-coherent scenarios it is assumed that the receiver
does not know the exact values of the entries of $H$ (other than
their statistical distribution). Another difference is that the
signal constellation $\V:=\{ \Phi_1,\Phi_2,\cdots, \Phi_L\}$ has
unitary constraints: $\Phi_k^* \Phi_k = I_M$ for
$k=1,2,\cdots,L$. The last equation simply states that the
columns of $\Phi_k$ form a ``unitary frame'', i.e., the column
vectors all have unit length in the complex vector space
$\mathbb{C}^T$ and the vectors are pairwise orthogonal.  The
scaled matrices $\sqrt{T} \Phi_k$, represent the codewords used
during the transmission.

The decoding task asks for the computation of the most likely
sent codeword $\hat{\Phi}$ given the received signal $R$. Under
the assumption of above model the ML decoder will
have to compute:
$$
\hat{\Phi}=\displaystyle \arg \max_{\Phi \in
  \{\Phi_1,\Phi_2,\cdots,\Phi_L\}} {\|R^*\Phi\|}
$$
for each received signal $R$ (see~\cite{ho00a}).

It has been shown in~\cite{ho00a} that the pairwise probability
of mistaking $\Phi_l$ for $\Phi_{l'}$ using ML decoding satisfies:
$$
P_{\Phi_l,\Phi_{l'}} \leq \frac{1}{2} \prod_{m=1}^M \left[1+
  \frac{(\rho T/M)^2(1-\delta_m^2 (\Phi_l^* \Phi_{l'}))}{4(1+\rho
    T/M)} \right]^{-N},
$$
where $\delta_m(\Phi_l^* \Phi_{l'})$ is the $m$-th singular value
of $\Phi_l^* \Phi_{l'}$. An important special case occurs when
$T=2M$. In this situation it is customary to represent the
unitary matrix $\Phi_k$ in the form:
\begin{equation}                         \label{specialform}
\Phi_k=\frac{\sqrt{2}}{2} \left(\begin{array}{c}
    I_M\\
    \Psi_k
  \end{array}\right).
\end{equation}
Note that by the definition of $\Phi_k$ the matrix $\Psi_k$ is a
$M \times M$ unitary matrix. In most of the literature mentioned
above researchers focus their attention on constellations having
the special form\eqr{specialform}. A differential modulation
scheme is discussed in~\cite{ho00} and this special form is used.
In this scheme, one does not send the identity matrix (the upper
part of the signal) every time. Instead of sending
$$\left(\begin{array}{cc}
    I_M\\
    \Psi_1\\
       \end{array}\right),   \left(\begin{array}{cc}
               I_M\\
               \Psi_2\\
       \end{array}\right), \left(\begin{array}{cc}
               I_M\\
               \Psi_3\\
       \end{array}\right), \cdots, $$
     one sends
     $$
     I_M, \Psi_1, \Psi_2\Psi_1,\Psi_3\Psi_2\Psi_1,\cdots.
     $$
     This increases the transmission rate by a factor of 2 to:
     $$
     \mathtt{R}=\frac{\log_2(L)}{M}=2\frac{\log_2(L)}{T}.
     $$
     Let $\Psi_{\tau}$ and $R_{\tau}$ denote the sent and
     received signals at time $\tau$ respectively, the ML decoder
     of the above differential space time modulation scheme will have to compute:
     $$
     \hat{\Psi}_{\tau}=\arg \min_{\Psi \in \V} \|R_{\tau}-\Psi
       R_{\tau-1} \|.
     $$
     For the constellations with the special
     form\eqr{specialform}, the pairwise error probability
     satisfies:
\begin{equation}   \label{pair-err}
P_{\Phi_l,\Phi_{l'}} \leq \frac{1}{2} \prod_{m=1}^M \left[1+
\frac{\rho^2\delta_m^2 (\Psi_l-\Psi_{l'})}{4(1+2\rho)}
\right]^{-N}.
\end{equation}
As explained in~\cite{ho00}, at high SNR scenarios the right side
of the above inequality is governed by the diversity product:
$$
\prod \V=\min_{l \neq l'} \frac{1}{2}
{|\det(\Psi_l-\Psi_{l'})|}^{1/M}.
$$
So the unitary differential modulation design criterion for a
non-coherent channel is to choose a constellation $\V$ such that
$\prod \V$ is as large as possible. We call a constellation $\V$
fully diverse if $\prod \V > 0$.

Interestingly enough the normalized complex orthogonal design
Alamouti codes and real orthogonal design codes can be used in
non-coherent scenarios too. Our work will be mainly about
non-coherent channels, however without any doubts the resulting
codes can be used for coherent channels too. In this paper we
will show how to construct space time codes from these schemes
using generalized phase shift keying (GPSK) signals. Our results
can be applied to generalized orthogonal design in~\cite{ja01}.

This paper is organized as follows: in Section~\ref{Section-2},
three series of $2$ dimensional GPSK constellations from complex
orthogonal designs will be presented. An algebraic calculation
will show that they have larger diversity products than the
original orthogonal design constellations. As a consequence of
larger diversity, the ML decoding of a GPSK constellation gives
better performance. In Section 3, we explicitly construct GPSK
constellations from the real orthogonal designs. Fast decoding
algorithms for the proposed constellations will be presented in
Section 2 and Section 3. Finally in Section 4 we will give the
conclusions and some future work.

\Section{GPSK Constellations from the Complex Orthogonal Design}
\label{Section-2}

A very simple yet interesting complex orthogonal scheme is
described in~\cite{al98}. We will consider the normalized version
of this proposed code, i.e., every element of this code is a
matrix given by
$$
\mathcal{O}(a,b)=\frac{1}{\sqrt{2}} \left(\begin{array}{cc}
    a & b \\
    -b^* & a^*
\end{array}\right),
$$
where $|a|^2=|b|^2=1$. Observe that $\mathcal{O}(a,b)$ is a
unitary matrix and a constellation $\mathcal{O}(n)$ with size
$L=n^2$ is obtained by letting $a$ and $b$ range over the set of
$n$-th roots of unity
$$
\vartheta=\{1,e^{2\pi i/n},\cdots,e^{2\pi i(n-1)/n}\},
$$
namely the entries of the matrix are chosen from scaled one
dimensional PSK signal set $\vartheta$. This would be the most
commonly implemented Alamouti's scheme. The diversity product of
the constellation $\mathcal{O}(n)$ is
$$
\prod \mathcal{O}(n)=\frac{\sqrt{2}}{2}\sin \frac{\pi}{n}.
$$
Note that $\mathcal{O}(n)$ is similar to $\vartheta$ in the sense
that all the elements in $\mathcal{O}(n)$ have unit energy and
every pair of elements differ only by the phases.
$\mathcal{O}(n)$ is a subset of the special unitary group
$$
SU(2)=\left\{\left(\begin{array}{cc}
      a&b\\
      -b^*&a^*\\
              \end{array}\right) \left|{|a|}^2+{|b|}^2=1 \right. \right\}.
$$
In this section we will present three series of
unitary constellations as finite subsets of $SU(2)$.

The basic principle to decode the constellation $\mathcal{O}(n)$
has been discussed in~\cite{te99}. We are going to describe this
decoding process in another way with more details and generalize
it to our constellations. Consider a non-coherent wireless
communication system modulated by $\mathcal{O}(n)$ with $2$
transmitting antennas and $N$ receiving antennas and assume the
differential unitary space time modulation~\cite{ho00} is used.
Let $X$ and $Y$ denote the received matrices at time block $\tau$
and $\tau+1$, respectively, then the ML decoder will perform the
following decoding task
$$
(\hat{a}, \hat{b})=\arg \min_{a, b \in \vartheta}
{\left\|Y-\frac{1}{\sqrt{2}}\left(\begin{array}{cc}
a&b\\
-b^*&a^*\\
\end{array}\right)X \right\|}^2.
$$
With simple matrix manipulations, one
can check that
$$
(\hat{a}, \hat{b})=\arg \min_{a, b \in
\vartheta} \sum_{i=1}^N
{\left\|\left(\begin{array}{c}
Y_{1i}\\
Y_{2i}\\
\end{array}\right)-\frac{1}{\sqrt{2}} \left(\begin{array}{cc}
a&b\\
-b^*&a^*\\
\end{array}\right) \left(\begin{array}{c}
            X_{1i}\\
            X_{2i}\\
    \end{array}\right) \right\|}^2
$$
$$
=\arg
\min_{a, b \in
\vartheta}
\sum_{i=1}^N
{\left\|\left(\begin{array}{c}
        Y_{1i}\\
        Y_{2i}^*\\
\end{array}\right)-\frac{1}{\sqrt{2}} \left(\begin{array}{cc}
X_{1i}&X_{2i}\\
X_{2i}^*&-X_{1i}^*\\
\end{array}\right) \left(\begin{array}{c}
            a\\
            b\\
    \end{array}\right) \right\|}^2
$$
$$
=\arg \min_{a, b \in \vartheta} \sum_{i=1}^N
\frac{2}{{|X_{1i}|}^2+{|X_{2i}|}^2} {\left\|\frac{1}{\sqrt{2}}
    \left(\begin{array}{cc}
        X_{1i}&X_{2i}\\
        X_{2i}^*&-X_{1i}^*\\
\end{array}\right)\left(\begin{array}{c}
Y_{1i}\\
Y_{2i}^*\\
\end{array}\right)- \frac{1}{2}({|X_{1i}|}^2+{|X_{2i}|}^2) \left(\begin{array}{c}
            a\\
            b\\
    \end{array}\right) \right\|}^2
$$
$$
=\arg \min_{a, b \in \vartheta} \sum_{i=1}^N
\frac{2}{{|X_{1i}|}^2+{|X_{2i}|}^2} \left({\left|
\frac{1}{\sqrt{2}}X_{1i}^*Y_{1i}+\frac{1}{\sqrt{2}}X_{2i}Y_{2i}^*-\frac{1}{2}({|X_{1i}|}^2+{|X_{2i}|}^2)a
\right|}^2+\right.
$$
$$
\left.{\left|
\frac{1}{\sqrt{2}}X_{2i}^*Y_{1i}+\frac{1}{\sqrt{2}}X_{1i}Y_{2i}^*-\frac{1}{2}({|X_{1i}|}^2+{|X_{2i}|}^2)b
\right|}^2\right).
$$
Further algebraic simplifications show that ML decoding is very
simple:
\begin{equation}  \label{MLD}
(\hat{a},\hat{b})=\arg \max_{a, b \in \vartheta} \left( Re(a
\sum_{i=1}^N \bar{Z}_i)+Re(b \sum_{i=1}^N \bar{W}_i) \right),
\end{equation}
where $Z_i=X_{1i}^*Y_{1i}+X_{2i}Y_{2i}^*$ and
$W_i=X_{2i}^*Y_{1i}+X_{1i}Y_{2i}^*$. Since $a$ and $b$ are
independent of each other ($a$ and $b$ can be chosen freely in
the set $\vartheta$), the evaluations above amount to
$$
\hat{a}=\arg \max_{a \in \vartheta} Re(a \sum_{i=1}^n \bar{Z}_i),
$$
and
$$
\hat{b}=\arg \max_{b \in \vartheta} Re(b \sum_{i=1}^n \bar{W}_i).
$$
Rewrite $a$ as $a=e^{2\pi j i/n}$ and $b$ as $b=e^{2\pi k i/n}$
and let
$$
\lfloor r \rceil=\lfloor r+1/2 \rfloor,
$$
i.e., $\lfloor r \rceil$ denotes the smaller of the closest
integers to $r$. The ML decoder will take the following simple
form:
\begin{equation}  \label{simple-MLD}
\hat{j}=\left\lfloor \frac{n \arg Z}{2\pi} \right\rceil, \qquad
\hat{k}=\left\lfloor \frac{n \arg W}{2\pi} \right\rceil,
\end{equation}
where $Z=\sum_{i=1}^n Z_i$ and $W=\sum_{i=1}^n W_i$.

Assume that a communication channel is modulated by one
dimensional PSK signals $\vartheta$. Let $a=e^{2\pi l/n} \in
\vartheta$ denote the sent signal at one time and $Z$ denote the
corresponding corrupted signal. The ML decoding will look for the
closest signal $\hat{a}=e^{2\pi \hat{j}/n}$ to $a$ in the signal
set $\vartheta$, i.e.,
\begin{equation}  \label{PSK-MLD}
\hat{j}=\left\lfloor \frac{n \arg Z}{2\pi} \right\rceil.
\end{equation}
Compare Formula\eqr{simple-MLD} and Formula\eqr{PSK-MLD}, one can
conclude that the decoding of space time code $\mathcal{O}(n)$ is
decomposable and can be implemented by taking one dimensional PSK
demodulation twice. The normalized Alamouti codes
$\mathcal{O}(n)$'s admit very simple decoding algorithms, since
roughly only $4N$ complex multiplications and $4N$ complex
additions are needed. We are going to introduce other complex
orthogonal codes using generalized PSK (GPSK) signals in the
sequel. Instead of allocating the same energy to every
transmitting antenna, we attempt to exploit the transmit
diversity and optimize the power allocation. Compared to
$\mathcal{O}(n)$, these codes have larger diversity, which
promise good performances with ML decoding. Due to the symmetry
from GPSK signals, these codes are also decomposable, therefore
fast decoding algorithms can be also applied to these codes. In
the following three subsections we present the constructions of
these three series of constellations.

\subsection{Construction 1}

Let $n$ be an even number and let $ 0 < r < \frac{\sqrt{2}}{2}$
be the root of the following equation
\begin{equation} \label{root}
\left(\frac{\sqrt{2}}{2}-r\right)^2+\left(\frac{\sqrt{2}}{2}-\sqrt{1-r^2}\right)^2=4r^2\sin^2\frac{2\pi}{n}.
\end{equation}
Consider the following sets of the scaled one dimensional PSK
signals:
$$
A_1(n)=\left\{\sqrt{2}/2e^{i\frac{2\pi}{n}k} | k=0,1,\cdots,n-1
\right\},
$$
$$
A_2(n)=\left\{re^{i\frac{4\pi}{n}k} |
  k=0,1,\cdots,\frac{n}{2}-1\right\},
$$
$$
A_3(n)=\left\{\sqrt{1-r^2}e^{i\frac{2\pi}{n}k} |
  k=0,1,\cdots,n-1\right\}.
$$
Consider the following subsets of $SU(2)$:
$$C_1(n)=\left\{\left(\begin{array}{cc}
      a & b \\
      -\bar{b} & \bar{a}
\end{array}\right)|a \in A_1(n) , b \in A_1(n) \right\},
$$
$$
C_2(n)=\left\{\left(\begin{array}{cc}
      a & b \\
      -\bar{b} & \bar{a}
\end{array}\right)|a \in A_2(n) , b \in A_3(n) \right\},
$$
$$
C_3(n)=\left\{\left(\begin{array}{cc}
      a & b \\
      -\bar{b} & \bar{a}
\end{array}\right)|a \in A_3(n) , b \in A_2(n) \right\}.
$$

\begin{thm}
  Let
  $$
  \V_1(n):=C_1(n) \bigcup C_2(n) \bigcup C_3(n),
  $$
  $\V_1(n)$ is a fully diverse constellation of $2n^2$
  elements with diversity product:
  $$
  \prod \V_1(n)=\min \left\{\frac{\sqrt{2}}{2} \sin
    \frac{\pi}{n},\frac{1}{2} \sqrt{(\frac{\sqrt{2}}{2}-r)^2+
      (\frac{\sqrt{2}}{2}-\sqrt{1-r^2})^2}\right\}.
  $$
\end{thm}

\begin{proof}

  Take two elements $A \in C_i$ and $B \in C_j$. Without loss of
  generality, we further assume $i \leq j$. If $ A , B \in
  C_1(n)$, then
  $$
  |\det(A-B)| \geq (\sqrt{2} \sin \frac{\pi}{n})^2=2 \sin^2
  \frac{\pi}{n}.
  $$
  If $A \in C_1(n), B \in C_2(n)$ or $ A \in C_1(n), B \in
  C_3(n)$ then
  $$|\det(A-B)| \geq (\frac
  {\sqrt{2}}{2}-r)^2+(\frac{\sqrt{2}}{2}-\sqrt{1-r^2})^2=2-\sqrt{2}r-\sqrt{2}\sqrt{1-r^2}.
  $$
  If $A \in C_2(n)$ and $B \in C_3(n)$ then
  $$
  |\det(A-B)| \geq 2(\sqrt{1-r^2}-r)^2.
  $$
  If $A,B \in C_3(n)$ or $A,B \in C_2(n)$ then
  $$
  |\det(A-B)| \geq (2r \sin \frac{2\pi}{n})^2.
  $$
  Use the fact that $r$ is the root of Equation\eqr{root} and
  compare the lower bounds of the diversity product in all the
  above cases, the claim in the theorem can be established.
\end{proof}

The following table shows how the diversity product of $\V_1(n)$
compares with the diversity product of $\mathcal{O}(n)$ when $n
\leq 12$. \vspace{0.3cm}

\begin{center}
\begin{tabular}{|c|c|c|c|c|} \hline
  $n$ &$r$& $\frac{\sqrt{2}}{2} \sin \frac{\pi}{n}=\prod \mathcal{O}(n)$ &
 $\frac{1}{2} \sqrt{(\frac{\sqrt{2}}{2}-r)^2+(\frac{\sqrt{2}}{2}-\sqrt{1-r^2})^2}$
  & $\prod \V_1(n)$ \\ \hline
  4 &0.259& 0.5 & 0.259 & 0.259 \\ \hline
  6 &0.284& 0.353 & 0.246 & 0.246 \\ \hline
  8 &0.321& 0.271 & 0.227 & 0.227 \\ \hline
  10 &0.360& 0.219 & 0.209 & 0.209 \\ \hline
  12 &0.386& 0.183 & 0.193 & 0.183 \\ \hline
\end{tabular}
\end{center}

\vspace{0.3cm}

\begin{co}
  For $n \geq 12 $,
  $$
  \prod \V_1(n) = \frac{\sqrt{2}}{2} \sin \frac{\pi}{n}.
  $$
\end{co}

\begin{proof}

  For $n \geq 12 $,
  $$
  \frac{1}{2}
  \sqrt{(\frac{\sqrt{2}}{2}-r)^2+(\frac{\sqrt{2}}{2}-\sqrt{1-r^2})^2}
  = r \sin \frac{2\pi}{n}=2r \sin \frac{\pi}{n} \cos\frac{\pi}{n}
  \geq 2 \times 0.386 \cos \frac{\pi}{12} \sin \frac{\pi}{n}
  $$
  $$
  \geq \frac{\sqrt{2}}{2} \sin \frac{\pi}{n}.
  $$
  Consequently we have
  $$
  \prod \V_1(n)=\min \left\{\frac{\sqrt{2}}{2} \sin
      \frac{\pi}{n},\frac{1}{2} \sqrt{(\frac{\sqrt{2}}{2}-r)^2+
        (\frac{\sqrt{2}}{2}-\sqrt{1-r^2})^2}\right\}=\frac{\sqrt{2}}{2} \sin \frac{\pi}{n}.
  $$

\end{proof}

The above corollary indicates that for $n \geq 12 $, the GPSK
constellation $\V_1(n)$ has the same diversity product as the
orthogonal constellation $\mathcal{O}(n)$, while it has twice as
many elements: $\V_1(n)$ has $2n^2$ elements whereas
$\mathcal{O}(n)$ has $n^2$ elements.

Similar to the case for $\mathcal{O}(n)$, the ML decoding of
$\V_1(n)$ boils down to Formula\eqr{MLD}. However we can not
separate the estimation of $a, b$ using the simple
Formula\eqr{simple-MLD}, because generally $a$ and $b$ are not
independent anymore. If we restrict the decoding evaluation in a
particular $C_i(n)$, then $\arg(a)$ and $\arg(b)$ can be chosen
freely. Namely within the restricted searching area $C_i(n)$, $a$
and $b$ are independent of each other. Assume the originally sent
codeword falls in $C_1(n)$, one can use Formula\eqr{simple-MLD}
to resolve a candidate $(\hat{a}_1, \hat{b}_1)$:
$$
\hat{j}=\left\lfloor \frac{n \arg Z}{2\pi} \right\rceil, \qquad
\hat{k}=\left\lfloor \frac{n \arg W}{2\pi} \right\rceil.
$$
Similarly we can have candidates $(\hat{a}_i, \hat{b}_i)$ for
$C_i(n)$, where $i=2,3$. The final ML decoder will resolve the
most likely sent codeword:
$$
(\hat{a}, \hat{b})=\arg \max_{\hat{a}_i, \hat{b}_i} \left
  (Re(a \bar{Z})+Re(b \bar{W}) \right).
$$
Again we conclude that the decoding of $\V_1(n)$ is decomposable
and therefore $\V_1(n)$ admits a simple decoding. The above
evaluations require roughly $4N$ complex multiplications and $4N$
complex additions, therefore the decoding of the GPSK
constellation $\V_1(n)$ has the same complexity as that of
$\mathcal{O}(n)$.

\subsection{Construction 2}

Let $n=2m$ and consider the following sets consisting of scaled
one dimensional PSK signals:
$$
A_1(n)=\left\{re^{i\frac{2\pi}{m}j} |
  j=0,1,\cdots,m-1\right\},
$$
$$
A_2(n)=\left\{\sqrt{1-r^2}e^{i(\frac{2\pi}{m}j+\frac{\pi}{m})}
  | j=0,1,\cdots,m-1\right\},
$$
$$
A_3(n)=\left\{\sqrt{1-r^2}e^{i\frac{2\pi}{m}j} |
  j=0,1,\cdots,m-1\right\},
$$
$$
A_4(n)=\left\{re^{i(\frac{2\pi}{m}j+\frac{\pi}{m})} |
  j=0,1,\cdots,m-1\right\},
$$
where
$$
r = \frac{1}{\sqrt{2\sin^2 \frac{\pi}{m}+2\sqrt{2}\sin
    \frac{\pi}{m}+2}}.
$$
Based on the above signal sets, construct the following
subsets of $SU(2)$:
$$
C_1(n)=\left\{\left(\begin{array}{cc}
      a & b \\
      -\bar{b} & \bar{a}
\end{array}\right)|a \in A_1(n),b \in A_2(n)\right\},
$$
$$
C_2(n)=\left\{\left(\begin{array}{cc}
      a & b \\
      -\bar{b} & \bar{a}
\end{array}\right)|a \in A_2(n),b \in A_1(n)\right\},
$$
$$
C_3(n)=\left\{\left(\begin{array}{cc}
      a & b \\
      -\bar{b} & \bar{a}
\end{array}\right)|a \in A_3(n),b \in A_4(n)\right\},
$$
$$
C_4(n)=\left\{\left(\begin{array}{cc}
      a & b \\
      -\bar{b} & \bar{a}
\end{array}\right)|a \in A_4(n),b \in A_3(n)\right\}.
$$

\begin{thm}
  Let
  $$
  \V_2(n):=C_1(n) \bigcup C_2(n) \bigcup C_3(n) \bigcup
  C_4(n),
  $$
  $\V_2(n)$ is a fully diverse constellation of $n^2$ elements
  with diversity product
  $$
  \prod \V_2(n)=\min\left\{r \sin \frac{2\pi}{n},\sin
    \frac{\pi}{n}\right\}.
  $$
\end{thm}

\begin{proof}

  Take two elements $A \in C_i(n)$ and $B \in C_j(n)$. Without
  loss of generality, we can further assume $i \leq j$. If $i=j$, we have
  $$
  |\det(A-B)| \geq
  |re^{i\frac{2\pi}{m}k}-re^{i\frac{2\pi}{m}(k+1)}|^2=4r^2 \sin^2
  \frac{\pi}{m}.
  $$
  If $A \in C_1(n),B \in C_2(n)$ or $A \in C_3(n),B \in
  C_4(n)$, one can check that
  $$
  |\det(A-B)| \geq
  2|re^{i\frac{2\pi}{m}k}-\sqrt{1-r^2}e^{i(\frac{2\pi}{m}k+\frac{\pi}{m})}|
  =2(1-2\sqrt{1-r^2}r\cos \frac{\pi}{m}).
  $$
  Similarly if $A \in C_1(n)$, $B \in C_3(n)$ or $A \in
  C_2(n)$, $B \in C_4(n)$,
  $$
  |\det(A-B)| \geq 2(\sqrt{1-r^2}-r)^2=2(1-2\sqrt{1-r^2}r).
  $$
  If $A \in C_1(n)$, $B \in C_4(n)$ or $A \in C_2(n)$, $B \in
  C_3(n)$,
  $$
  |\det(A-B)| \geq |re^{i\frac{2\pi}{m}k}-re^{i\frac{2\pi
      k}{m}+\frac{\pi}{m}}|^2+|
  \sqrt{1-r^2}e^{i\frac{2\pi}{m}k}-\sqrt{1-r^2}e^{i\frac{2\pi
      k}{m}+\frac{\pi}{m}}|^2 =4\sin^2 \frac{\pi}{2m}.
  $$
  It follows from the definition of $r$ that
  $$
  2(1-2\sqrt{1-r^2}r)=4r^2\sin^2 \frac{\pi}{m},
  $$
  and naturally we will have
  $$
  2(1-2\sqrt{1-r^2}r\cos \frac{\pi}{m}) \ge
  2(1-2\sqrt{1-r^2}r).
  $$
  Compare the lower bounds in all the cases and take the
  minimum of them, we establish the claim in the theorem.
\end{proof}

The following table shows how the diversity product of $\V_2(n)$
compares with that of $\mathcal{O}(n)$ when $n \leq 14$.

\vspace{0.3cm}

$$\begin{array}{|c|c|c|c|c|c|} \hline n & r & r \sin \frac{2\pi}{n} & \sin \frac{\pi}{n} & \prod \V_2(n) &
  \prod \mathcal{O}(n) \\ \hline 4 & 0.383 & 0.383 & 0.707 & 0.383 & 0.5\\
  \hline 6 & 0.410 & 0.355 & 0.500 & 0.355 & 0.354 \\ \hline 8 &
  0.447 & 0.316 & 0.355 & 0.316 & 0.271 \\ \hline 10 & 0.479 &
  0.282 & 0.309 & 0.282 & 0.219 \\ \hline 12 & 0.505 & 0.253 &
  0.259 & 0.253 & 0.183 \\ \hline 14 & 0.527 & 0.229 & 0.222 &
  0.222 & 0.157 \\ \hline
\end{array}
$$

\vspace{0.3cm}

\begin{co}
  For $n \geq 14$,
  $$\prod\V_2(n)=\sin \frac{\pi}{n}.$$
\end{co}
\begin{proof}
  According to the definition of $r$,
  $$
  2(1-2\sqrt{1-r^2}r)=4r^2\sin^2 \frac{\pi}{m}.
  $$
  One can also check that
  $$
  4r^2 \sin^2 \frac{\pi}{m}=16r^2 \sin^2 \frac{\pi}{2m} \cos^2
  \frac{\pi}{2m}=4r^2 \cos^2 \frac{\pi}{2m} (4 \sin^2
  \frac{\pi}{2m})
  $$
  $$
  =\frac{4 \cos^2 \frac{\pi}{2m}}{2 \sin^2
    \frac{\pi}{m}+2\sqrt{2}\sin \frac{\pi}{m}+2}(4 \sin^2
  \frac{\pi}{2m}) \geq \frac{4 \cos^2 \frac{\pi}{14}}{2 \sin^2
    \frac{\pi}{7}+2\sqrt{2}\sin \frac{\pi}{7}+2} (4 \sin^2
  \frac{\pi}{2m}) \ge 4 \sin^2 \frac{\pi}{2m}.
  $$
  Consequently one concludes that
  $$
  \prod \V_2(n)=\min\left\{r \sin \frac{\pi}{m},\sin
    \frac{\pi}{2m} \right\}=\sin \frac{\pi}{n}.
  $$
\end{proof}

For $n \geq 14$, the constellation $\V_2(n)$ has as many elements
as the constellation $\mathcal{O}(n)$, however its diversity
prodcut is larger than that of $\mathcal{O}(n)$ by a factor of
$\sqrt{2}$, i.e., $\prod \V_2(n)= \sqrt{2} \prod \mathcal{O}(n)$.
The corollary indicates that we could allocate different power to
the transmitting antennas to achieve more reliable transmission,
while still keeping the total energy. Similar to the case for
$\V_1(n)$, one can apply exactly the same algorithm to achieve
the ML decoding for $\V_2(n)$. It can be easily seen that these
two algorithms have the same complexity.

\subsection{Construction 3}

We will take further efforts to explore the subsets of $SU(2)$.
In the following we are going to describe a series of unitary
constellation $\V_3(n)$. For the sample program implemented to
construct $\V_3(n)$, we refer to~\cite{ha03u2}. Now for given
integers $n > 0$ and $ 0 \leq k \leq n$, we define
$$
N_0=1,N_k=\frac{\pi}{\left\lfloor \arcsin \frac{\sin
      \frac{\pi}{4n}}{\cos \frac{(n-k)\pi}{2n}} \right\rfloor},
k=1,2,\cdots,n,
$$
$$
M_k=\frac{\pi}{\left\lfloor \arcsin \frac{\sin
      \frac{\pi}{4n}}{\sin \frac{(n-k)\pi}{2n}} \right\rfloor},
k=0,1,\cdots,n-1, M_n=1,
$$
$$
C_k(n)=\left\{\left(\begin{array}{cc}
      a_{k,j}&b_{k,l}\\
      -\bar{b}_{k,l}&\bar{a}_{k,j}\\
                \end{array}\right)|a_{k,j}=\cos
              \frac{(n-k)\pi}{2n} e^{i \frac{2j\pi}{N_k}},
              b_{k,l}=\sin \frac{(n-k)\pi}{2n} e^{i
                \frac{2l\pi}{M_k}} \right\}.
            $$

\begin{thm}   \label{construction}
  Let
  $$
  \V_3(n)=\bigcup_{k=0}^n C_k(n),
  $$
  $\V_3(n)$ is a fully diverse constellation of
  $\displaystyle\sum_{k=0}^n M_kN_k $ elements with diversity
  product
  $$
  \prod \V_3(n) = \sin \frac{\pi}{4n}.
  $$

\end{thm}

\begin{proof}

  Pick two distinct elements $A, B \in \V$,
  $$
  A=\left(\begin{array}{cc}
      a&b\\
      -\bar{b}&\bar{a}
          \end{array}\right), \qquad \qquad \qquad \qquad B=\left(\begin{array}{cc}
             c&d\\
             -\bar{d}&\bar{c}
          \end{array}\right).
        $$
        One can verify that
        $$
        |\det(A-B)|=\det(A-B)=|a-c|^2+|b-d|^2.
        $$
        So if $|a| \ne |c|$, then we have
        $$
        |\det(A-B)|=|a-c|^2+|b-d|^2 \geq
        (|a|-|c|)^2+(|b|-|d|)^2 \geq 2-2\cos \frac{\pi}{2n},
        $$
        and it can be verified that the equality holds if
        there is a $k \in \{0, 1, \cdots, n\}$ such that $A \in
        C_k(n)$ and $B \in C_{k+1}(n)$ or alternatively $B \in
        C_k(n)$ and $A \in C_{k+1}(n)$. In the case that $|a| =
        |c|$, we will have
        $$
        |\det(A-B)|=|a-c|^2+|b-d|^2 \geq \max \{|a-c|^2,|b-d|^2\} \geq 2-2\cos
        \frac{\pi}{2n}.
        $$
        Therefore for all the cases, we have
        $$
        |\det(A-B)| \geq 2-2 \cos \frac{\pi}{2n}.
        $$
        One checks that the lower bound for each case can
        be reached.  So for constellation $\V_3(n)$, it follows
        that
        $$
        \prod \V_3(n)=\frac{1}{2} (2-2 \cos
        \frac{\pi}{2n})^{\frac{1}{2}}=\sin \frac{\pi}{4n}.
        $$

\end{proof}

\begin{co}  \label{V3-asymptotic}
  When $n \rightarrow \infty$, $\V_3(n)$ will have have $O(n^3)$
  elements and have the diversity product $O(\frac{1}{n})$.
\end{co}

\begin{proof}
  For $k=1, 2, \cdots, n$,
  $$
  N_k=\frac{\pi}{\left\lfloor \arcsin \frac{\sin
        \frac{\pi}{4n}}{\cos \frac{(n-k)\pi}{2n}} \right\rfloor}
  \leq \frac{\pi}{\left\lfloor \arcsin \sin \frac{\pi}{4n}
    \right\rfloor} \leq 4n.
  $$
  Similarly for $k=0, 1, \cdots, n-1$,
  $$
  M_k=\frac{\pi}{\left\lfloor \arcsin \frac{\sin
        \frac{\pi}{4n}}{\sin \frac{(n-k)\pi}{2n}} \right\rfloor}
  \leq \frac{\pi}{\left\lfloor \arcsin \sin \frac{\pi}{4n}
    \right\rfloor} \leq 4n.
  $$
  Hence we derive an asymptotic upper bound for the
  cardinality of $\V_3(n)$,
  $$
  |\V_3(n)|=\sum_{k=0}^n|C_k(n)|=\sum_{k=0}^n M_k N_k \leq
  2n+16(n-1)n^2 \leq O(n^3).
  $$
  In the following we are going to derive an asymptotic lower
  bound for $|\V_3(n)|$. First pick two real numbers $\alpha,
  \beta$ such that
  $$
  0 < \alpha < \beta <1.
  $$
  For $k$ such that $ \alpha n < k < \beta n$ (such $k$ always
  exists when $n$ is large enough), we have
  $$
  N_k=\frac{\pi}{\left\lfloor \arcsin \frac{\sin
        \frac{\pi}{4n}}{\cos \frac{(n-k)\pi}{2n}} \right\rfloor}
  \geq \frac{\pi}{\left\lfloor \arcsin \frac{\sin
        \frac{\pi}{4n}}{\cos \frac{(1-\alpha)\pi}{2}}
    \right\rfloor}.
  $$
  Utilizing the fact that
  $$
  \sin x \sim x \sim \arcsin x,
  $$
  when $x$ is close to $0$. When $n$ is large enough, one derives
  $$
  N_k \geq 4\cos \frac{(1-\alpha)\pi}{2}n.
  $$
  Similarly
  $$
  M_k \geq 4\sin \frac{(1-\beta)\pi}{2}n.
  $$
  Therefore we have an asymptotic lower bound for the
  cardinality of $\V_3(n)$,
  $$
  |\V_3(n)| \geq \sum_{\alpha n < k < \beta n} M_k*N_k \geq 16
  (\beta-\alpha) \sin \frac{(1-\beta)\pi}{2} \cos
  \frac{(1-\alpha)\pi}{2} n^3 \geq O(n^3).
  $$
  Combining the upper bound and lower bound, we conclude that
  $\V_3(n)$ has $O(n^3)$ elements with diversity product
  $\sin (\pi/(4n)) = O(1/n)$.
\end{proof}

This corollary indicates that asymptotically $\V_3(n)$ will have
much better diversity product compared to $\mathcal{O}(n),
\V_1(n), \V_2(n)$, because the other three constellations
asymptotically will have $O(n^2)$ elements and have diversity
product $O(\frac{1}{n})$. Observe that $SU(2)$ in fact can be
viewed as $3$ dimensional unit sphere. Finding a constellation
with the optimal diversity product can be converted to a sphere
packing problem on this unit sphere. For $n$ points on a $3$
dimensional unit sphere, the largest minimum distance one can
hope for is asymptotically $O(1/n^{1/3})$. Therefore when $n$
becomes large, asymptotically $\V_3(n)$ is the best constellation
over all the subsets of $SU(2)$. Table~\ref{V3} shows the
diversity product of constellation $\V_3(n)$.

\begin{table}
\begin{center}
\begin{tabular}{|c|c|c|}
  \hline
  $n$ &size & diversity product \\ \hline
  3&124 & 0.259 \\ \hline
  4&293 & 0.195 \\ \hline
  5&582 & 0.156 \\ \hline
  6&974 & 0.131 \\ \hline
  7&1640 & 0.112 \\ \hline
  8&2438 & 0.098 \\ \hline
  9&3510 & 0.087 \\ \hline
  10&4898& 0.078 \\ \hline
  11&6516 & 0.071 \\ \hline
  12&8433& 0.065 \\ \hline
  13&10770 & 0.060 \\ \hline
\end{tabular}
\end{center}
\caption{the diversity product of $\V_3(n)$} \label{V3}
\end{table}

The ML decoding of $\V_3(n)$ is also decomposable and it roughly
requires $4N+2n$ complex multiplications and $4N+2n$ complex
additions, therefore asymptotically it requires $O(N)+O(L^{1/3})$
complex multiplications and additions.  One can see from
Table~\ref{V3} that for reasonably small values of $n$ the
constellation size grows rapidly. For example when the
constellation size $L$ is already huge (10770), $n$ is still
reasonably small (13). So even if we are dealing with a high rate
constellation, the decoding process of $\V_3(n)$ is still very
simple.

\begin{figure}[ht]
\centerline{\psfig{figure=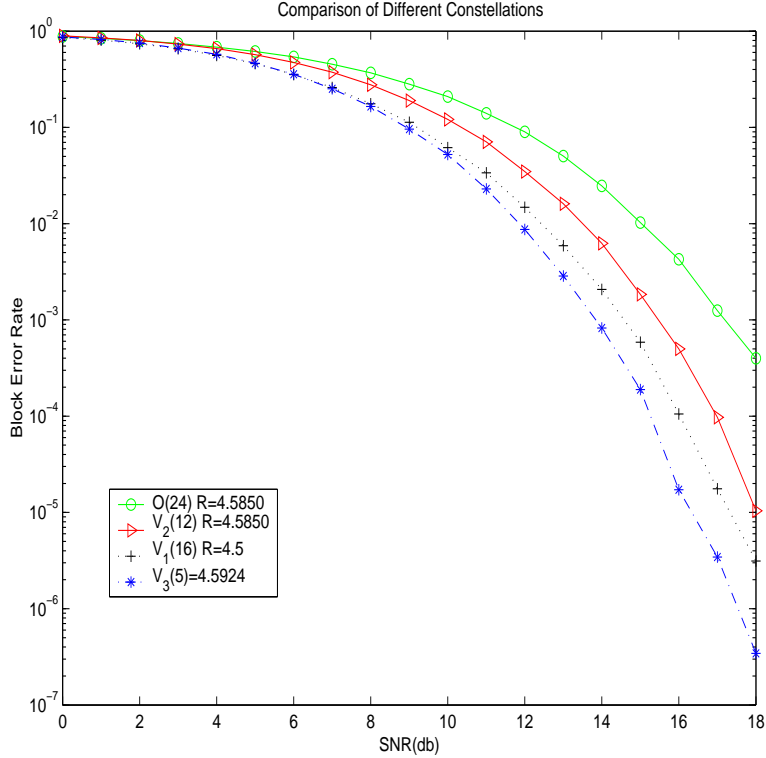,width=4in,height=4in}}
\caption{Performance of different codes from GPSK signals, $M$=2,
$N$=12} \label{figure-one}
\end{figure}

We compare different codes from GPSK signals with transmission
rate around $4.5$ in Figure~\ref{figure-one}. One can see that
$\V_3(5)$ outperforms other constellations even with the highest
transmission rate. Of course the decoding of $\V_3(5)$ is a little
more complex than other constellations, however the sacrifice in
decoding efficiency is worthwhile for the remarkably gained
performance.

\begin{figure}[ht]
  \centerline{\psfig{figure=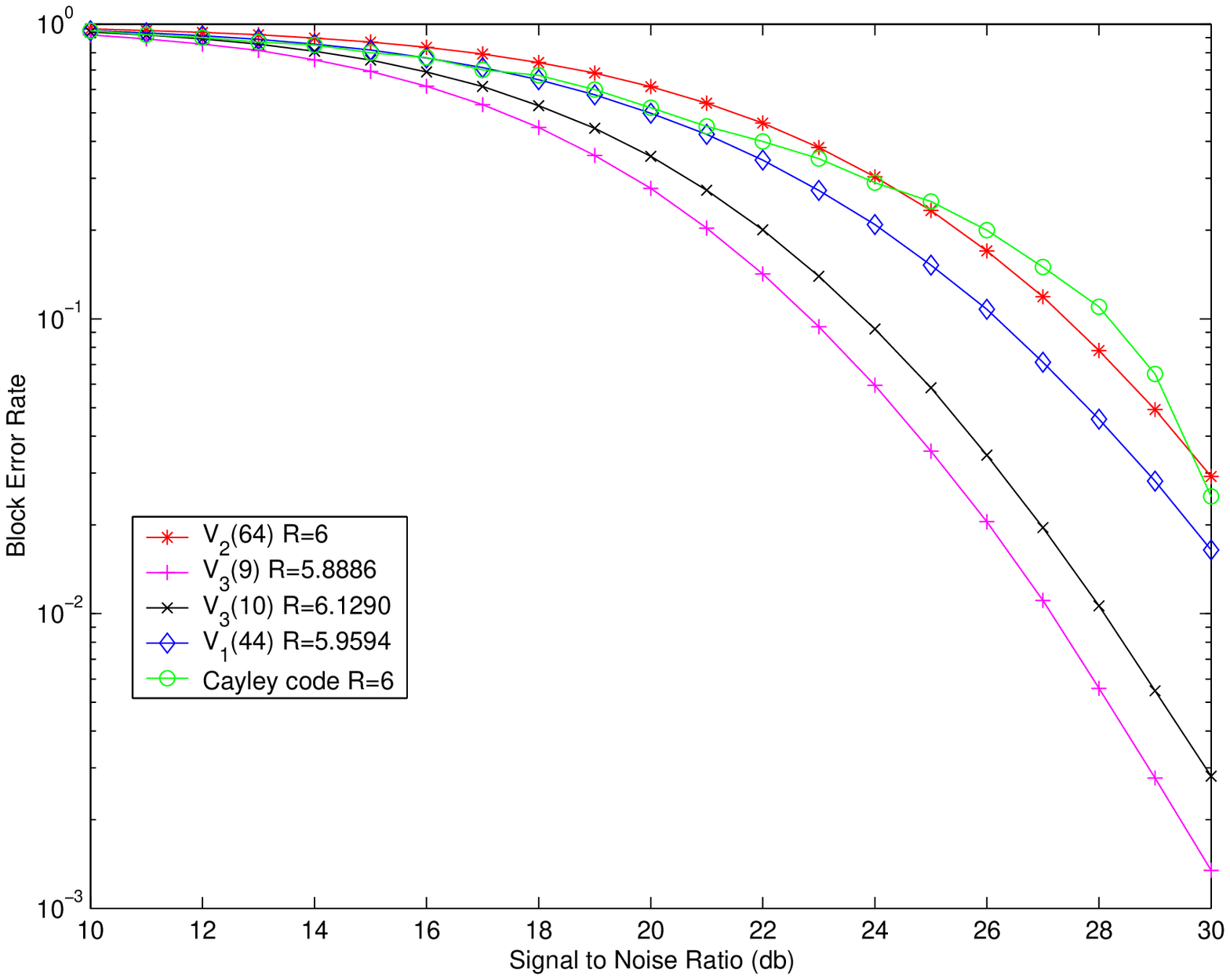,width=4in,height=4in}}
\caption{Performances of GPSK constellations and a Cayley code,
$M$=2, $N$=2} \label{figure-two}
\end{figure}

We compare the GPSK constellations with the Cayley
codes~\cite{ha02a} in Figure~\ref{figure-two}. One can see that
$\V_1(44)$ with transmission rate $R=5.9594$ has already a gain
of about $2$db compared to the Cayley code, and the performances
of $\V_3(9) \; \mbox{with}\; R=5.8886$ and $\V_3(10)
\;\mbox{with}\; R=6.1290$ are even more remarkable. Note that all
the GPSK constellations admit the presented simple decoding
alrogithms. $\V_1(44)$ or $\V_2(64)$ only needs about $8$ complex
multiplications and $8$ complex additions to decode one codeword,
$\V_3(9)$ or $\V_3(10)$ needs about $28$ complex multiplications
and $28$ complex additions to decode.

\Section{GPSK Constellations from the Real Orthogonal Designs}

In~\cite{ta99} space-time block coding has been discussed for
coherent channels. The authors classify all the real orthogonal
designs and show that orthogonal design codes only exist for $2,
4$ or $8$ dimension. We show that the normalized version of the
real orthogonal codes can be used in unitary space time
modulation for non-coherent scenarios too. We are going to
present explicit implementations of the real orthogonal designs,
which feature simple decoding algorithms. For simplicity we
describe our idea using a $4$ dimensional real orthogonal design:
\begin{equation}  \label{real4}
S=\left(\begin{array}{cccc}
    s_1&s_2&s_3&s_4\\
    -s_2&s_1&-s_4&s_3\\
    -s_3&s_4&s_1&-s_2\\
    -s_4&-s_3&s_2&s_1\\
      \end{array} \right).
\end{equation}
    One natural idea is to transform a $2$ dimensional complex
    orthogonal code $\mathcal{O}(n)$ to a $4$ dimensional real orthogonal
    code.  Recall that $\mathcal{O}(n)$ is a set of matrices
    $$
    \mathcal{O}(a,b)=\frac{1}{\sqrt{2}} \left(\begin{array}{cc}
        a & b \\
        -b^* & a^*
\end{array}\right)
$$
with $a, b$ ranging over $\vartheta=\{1,e^{2\pi
  i/n},\cdots,e^{2\pi i(n-1)/n}\}$. Let $s_1=Re(a), s_2=Im(a),
s_3=Re(b), s_4=Im(b)$, we will have a $4$ dimensional real
orthogonal code with the same diversity product as
$\mathcal{O}(n)$.

This code also features a simple decoding algorithm. The simple
decoding of real orthogonal design codes has been explained
in~\cite{ta99}, in the following we will discuss the decoding
process in more details such that it can be generalized to our
constellations straightforwardly. First we introduce some
notations to be used. The superscript $T$ will be the notation
for a matrix transpose. For two vectors $F$ and $G$, let $F \cdot
G$ denote the dot product of $F$ and $G$. For a vector $F$, $F_j$
denotes the $j$-th element of $F$. For a matrix $F$, $F_j$ denotes
the $j$-th column vector. Let $B$ be the mapping
$$
B(x_1, x_2, x_3, x_4)=\left( \begin{array}{cccc}
       x_1&x_2&x_3&x_4\\
       x_2&-x_1&x_4&-x_3\\
       x_3&-x_4&-x_1&x_2\\
       x_4&x_3&-x_2&-x_1\\
       \end{array} \right),
$$
immediately one checks that
\begin{equation}  \label{my-equation}
S \left( \begin{array}{c}
          x_1\\
          x_2\\
          x_3\\
          x_4\\
  \end{array} \right)= B(x_1, x_2, x_3, x_4) \left( \begin{array}{c}
          s_1\\
          s_2\\
          s_3\\
          s_4\\
  \end{array} \right).
\end{equation}

Assume that the differential unitary space time modulation is
used for a wireless communication system with $M$ transmitting
antennas and $N$ receiving antennas. Let $X, Y$ denote the
received $M \times N$ matrices at time block $\tau-1$ and $\tau$,
respectively. The ML decoder will make the following estimation:
$$
\hat{S}=\arg \min_S {\|Y-SX\|}^2.
$$
We separate $X, Y$ into the real parts and the imaginary parts
with
$$
X=F+iG, \qquad Y=P+iQ,
$$
then we have
$$
\hat{S}=\arg \min_S \left( {\|P-SF\|}^2+{\|Q-SG\|}^2
\right)= \arg \min_S \left( \sum_{j=1}^N
  {\|P_j-SF_j\|}^2+\sum_{j=1}^N {\|Q_j-SG_j\|}^2 \right).
$$
Utilize Equation\eqr{my-equation}, we have
$$
\hat{S}=\arg \min_{s_1, s_2, s_3, s_4} \left( \sum_{j=1}^N
{\left\|P_j-B(F_j) \left(\begin{array}{c}
                                                         s_1\\
                                                         s_2\\
                                                         s_3\\
                                                         s_4\\
                                              \end{array}
                                           \right)
                                         \right\|}^2+\sum_{j=1}^N {\left\|Q_j-B(G_j) \left( \begin{array}{c}
                                                         s_1\\
                                                         s_2\\
                                                         s_3\\
                                                         s_4\\
                                              \end{array} \right)
                                          \right\|}^2 \right).
$$
Since $B(F_j), B(G_j)$ are orthogonal matrices, simple
algebraic manipulations will lead the above evaluation to
$$
\arg \min_{s_1, s_2, s_3, s_4} \left(\sum_{j=1}^N
\frac{1}{{|F_j|}^2}
  {\left\|B^T(F_j)P_j-{|F_j|}^2 \left(\begin{array}{c}
                                                         s_1\\
                                                         s_2\\
                                                         s_3\\
                                                         s_4\\
                                              \end{array}\right) \right\|}^2+
  \sum_{j=1}^N \frac{1}{{|G_j|}^2} {\left\|B^T(G_j)Q_j-{|G_j|}^2 \left(\begin{array}{c}
                                                         s_1\\
                                                         s_2\\
                                                         s_3\\
                                                         s_4\\
                                              \end{array} \right)\right\|}^2
                                              \right)
$$
$$
=\arg \max_{s_1, s_2, s_3, s_4} (s_1, s_2, s_3, s_4) \cdot (\sum_{j=1}^N B^T(F_j)
P_j+ \sum_{j=1}^N B^T(G_j) Q_j)=\arg \max_{s_j} \sum_{j=1}^N s_jU_j,
$$
where $U_j=\sum_{j=1}^N B^T(F_j) P_j+ \sum_{j=1}^N B^T(G_j)
Q_j$.  Since the construction of our code is based on
$\mathcal{O}(n)$ as described above, one can check that $(s_1,
s_2)$ in fact is independent of $(s_3, s_4)$. Rewrite $s_1+i
s_2=\frac{1}{\sqrt{2}}e^{2k\pi i/n}$ and $s_3+i
s_4=\frac{1}{\sqrt{2}}e^{2l\pi i/n}$, then from
$$
\hat{S}= \arg \max_{s_1, s_2, s_3, s_4} \left( (s_1, s_2)
  \cdot (U_1, U_2)^T+(s_3, s_4) \cdot (U_3, U_4)^T \right),
$$
we conclude that the ML decoding of this GPSK constellation is
decomposable and can be boiled down to the following simple form:
$$
\hat{k}=\left\lfloor \frac{2\pi \arg (U_1+i U_2)}{n}
\right\rceil, \qquad \hat{l}=\left\lfloor \frac{2\pi \arg (U_3+i
U_4)}{n} \right\rceil.
$$
Generally speaking the use of the proposed codes above for a
wireless communication system with $M$ transmitting antennas and
$N$ receiving antennas will take $8M^2N$ real multiplications and
$8M^2N$ real additions to decode one codeword, which is very
simple.

Apparently we can apply the same idea for $\V_i(n)$ ($i=1,2,3$)
to construct GPSK real orthogonal constellations. The ML decoder
of the corresponding codes will take similar approaches as in the
$\mathcal{O}(n)$ case, except that one should notice that $(s_1,
s_2)$ is not independent oof $(s_3, s_4)$ anymore. In this case
one can restrict the searching area to be the subsets $C_i(n)$
``locally'' and apply the similar techniques as we did in the
complex orthogonal design case to achieve the ML decoding, then
the corresponding codes will admit simple decoding algorithms
too. Also we conclude that codes from $\V_1(n)$ or $\V_2(n)$ are
of the same complexity as codes from $\mathcal{O}(n)$. For the
codes from $\V_3(n)$, the complexity is increased for high
transmission rates, however they will have more pronounced
performances.

In the sequel we are going to present a series of $8$ dimensional
GPSK real orthogonal constellation $\V_4(n)$.  Consider the $8$
dimensional orthogonal design:
\begin{equation}   \label{eight-design}
\left( \begin{array}{cccccccc}
          s_1&s_2&s_3&s_4&s_5&s_6&s_7&s_8\\
          -s_2&s_1&s_4&-s_3&s_6&-s_5&-s_8&s_7\\
          -s_3&-s_4&s_1&s_2&s_7&s_8&-s_5&-s_6\\
          -s_4&s_3&-s_2&s_1&s_8&-s_7&s_6&-s_5\\
          -s_5&-s_6&-s_7&-s_8&s_1&s_2&s_3&s_4\\
          -s_6&s_5&-s_8&s_7&-s_2&s_1&-s_4&s_3\\
          -s_7&s_8&s_5&-s_6&-s_3&s_4&s_1&-s_2\\
          -s_8&-s_7&s_6&s_5&-s_4&-s_3&s_2&s_1\\
         \end{array} \right).
\end{equation}
Similar to the $4$ dimensional real orthogonal constellation
case, one can obtain $8$ dimensional unitary codes by
transforming $V(n) \times W(n)$, where $V$ or $W$ denotes any one
of $\mathcal{O}, \V_1, \V_2, \V_3$ and $\times$ denotes the
Cartesian product. For this implementation, $s_i$ can be assigned
to be the scaled version of the real or imaginary part of ``$a$''
or ``$b$'' as in $4$ dimensional case. However as in the complex
orthogonal case $\V_3(n)$ motivates us to explore more densely
packed constellations.

The problem of constructing an $8$ dimensional real orthogonal
code with the maximal diversity product is equivalent to the
packing problem on a $7$ dimensional unit sphere. Therefore any
currently existing results in packing problem on a $7$
dimensional unit sphere can be ``borrowed'' for the real
orthogonal code construction. However decoding such codes using
exhaustive search will be unpractical for high transmission
rates. In the sequel we are going to present a series of $8$
dimensional GPSK orthogonal codes $\V_4(n)$ featuring simple
decoding algorithms.

The basic idea is that instead of considering the problem of
packing the real vector $\mathbf{s}=(s_1, s_2, \cdots, s_8)$ on a
unit sphere, we consider the problem of packing the complex vector
$\mathbf{z}=(z_1, z_2, z_3, z_4)$ with constraints $\sum_{j=1}^4
{|z_j|}^2=1$, where $z_j=s_{2j-1}+i s_{2j}$. Since the distance of
two complex vectors satisfies
$$
\sum_{j=1}^4 {|z_j-w_j|}^2 \geq \sum_{j=1}^4
{(|z_j|-|w_j|)}^2,
$$
we can pack the amplitude vector $(|z_1|, |z_2|, |z_3|, |z_4|)$
first, then pack the argument vector
$(\arg(z_1),\arg(z_2),\arg(z_3),\arg(z_4))$. To explain this
implementation in more details, we represent $z$ using the polar
coordinates:
$$
\begin{tabular}{l}
  $z_1=\cos \theta_1 e^{i \gamma_1}$, \\
  $z_2=\sin \theta_1 \cos \theta_2 e^{i \gamma_2}$,\\
  $z_3=\sin \theta_1 \sin \theta_2 \cos \theta_3 e^{i \gamma_3}$,\\
  $z_4=\sin \theta_1 \sin \theta_2 \sin \theta_3 e^{i
  \gamma_4}$.\\
\end{tabular}
$$
Let $\theta_j$ run over $m_j+1$ evenly distributed discrete
values from $0$ to $\pi/2$:
$$
\theta_j \in \{k \frac{\pi}{2m_j}|k=0, 1, \cdots, m_j\},
$$
and let $\gamma_j$ run over $n_j$ evenly distributed discrete
values from $0$ to $2\pi$:
$$
\gamma_j \in \{k \frac{2\pi}{n_j}|k=0, 1, \cdots, n_j-1\}.
$$
Now we have a finite set of complex vectors.

Let $\mathbf{v}=(v_1, v_2, v_3, v_4)$ and $\mathbf{w}=(w_1, w_2,
w_3, w_4)$ be two distinct resulting complex vectors. If $|v_1|
\neq |w_1|$, one can check that
$$
{|\mathbf{v}-\mathbf{w}|}^2 \geq 4\sin^2 (\pi/(4m_1)).
$$
If $|v_1| = |w_1|$ but $|v_2| \neq |w_2|$, we will have
$$
{|\mathbf{v}-\mathbf{w}|}^2 \geq 4(1-{|v_1|}^2)\sin^2
(\pi/(4m_2)).
$$
In the case that $|v_1| = |w_1|$ and $|v_2| = |w_2|$ but $|v_3|
\neq |w_3|$, similar algebraic calculations will give a lower
bound
$$
{|\mathbf{v}-\mathbf{w}|}^2 \geq 4(1-{|v_1|}^2-{|v_2|}^2)\sin^2
(\pi/(4m_3)).
$$
We further consider the case $|v_j|=|w_j|$ for $j=1, 2, 3, 4$.
In this case, if $v_j \neq w_j$ for some $j$, we can have
$$
{|\mathbf{v}-\mathbf{w}|}^2 \geq 4{|v_j|}^2 \sin^2 (\pi/{n_j}).
$$
It is easy to check that the lower bound can be reached for some
special pair of $z, w$ in all the cases. Based on the $8$
dimensional real orthogonal design~\eqr{eight-design}, assign
$s_{2j-1}=Re(z_j)$ and $s_{2j}=Im(z_j)$ for $j=1,2,3,4$. Now we
have a finite unitary constellation whose diversity product is
$$
\min_{\mathbf{v}, \mathbf{w}} \sqrt{|\mathbf{v}-\mathbf{w}|}/2.
$$
\begin{figure}[ht]
  \centerline{\psfig{figure=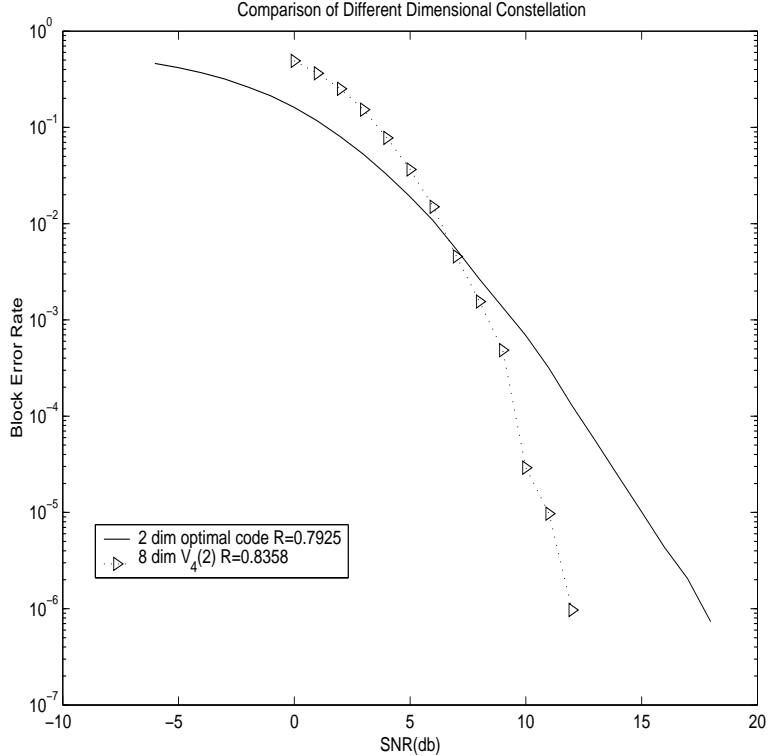,width=4in,height=4in}}
\caption{Performances of a GPSK constellation and a code with
  optimal diversity product, $N$=2}
\label{figure-three}
\end{figure}

However if we choose $m_j, n_j$ arbitrarily, often we will end up
having a constellation with a small diversity. In the sequel we
are going to add suitable constraints to the choices of $m_j$ and
$n_j$ to guarantee that we have a large diversity product. We
present the algorithm step by step. For the sample program of
this algorithm, we refer to~\cite{ha03u2}.

\vspace{0.3cm}

\textbf{Algorithm to Construct $\V_4(n)$:} (given the input $n$)
\begin{enumerate}
\item Fix $m_1=n$ and let $\theta_1$ run over $n+1$ even
  distributed discrete values from $0$ to $\pi/2$ (consequently
  $|z_1|$ runs over $n+1$ discrete values from $0$ to $1$).
\item For any fixed $|z_1|$, take $m_2$ to be the largest integer
  such that
  $$
  4(1-{|z_1|}^2)\sin^2 (\pi/(4m_2)) \geq 4\sin^2 (\pi/(4n)).
  $$
  Let $\theta_2$ run over $m_2+1$ evenly distributed discrete
  values from $0$ to $\pi/2$.
\item For any fixed $|z_1|, |z_2|$, take $m_3$ to be the largest
  integer such that
  $$
  4(1-{|z_1|}^2-{|z_2|}^2)\sin^2 (\pi/(4m_3)) \geq 4\sin^2
  (\pi/(4n)).
  $$
  Let $\theta_3$ run over $m_3+1$ evenly distributed discrete
  values from $0$ to $\pi/2$.
\item For any fixed $|z_1|, |z_2|, |z_3|$, take $n_j$ to be the
  largest integer such that
  $$
  4{|z_j|}^2 \sin^2 (\pi/{n_j}) \geq 4\sin^2 (\pi/(4n)).
  $$
  Let $\gamma_j$ run over $n_j$ evenly distributed discrete
  values from $0$ to $2 \pi$.
\item The steps above result in a finite set of complex vectors
  $\mathbf{z}=(z_1, z_2, z_3, z_4)$. Based on the $8$ dimensional real
  orthogonal design~\eqr{eight-design}, assign $s_{2j-1}=Re(z_j)$
  and $s_{2j}=Im(z_j)$ with $j=1,2,3,4$. Now we have a fully
  diverse unitary constellation with diversity product $\sin
  (\pi/ (4n))$.
\end{enumerate}

We call this constellation $\V_4(n)$. Recall that for the packing
problem with $n$ points on a $7$ dimensional unit sphere,
asymptotically the largest minimum distance one can hope for is
$O(1/n^7)$. The following theorem indicates that asymptotically
$\V_4(n)$ are the best constellations from $8$ dimensional real
orthogonal designs. Although the proof of this theorem is quite
tedious, the basic idea is very similar to
Corollary~\ref{V3-asymptotic}, thus we skip the proof.
\begin{thm}
  $\V_4(n)$ is a fully diverse constellation with diversity product $\sin
  (\pi/ (4n))$. When $n \rightarrow
  \infty$, $\V_4(n)$ will have have $O(n^7)$ elements and have
  the diversity product $O(1/n)$.
\end{thm}

Applying the same analysis for the constellation $\V_3(n)$, we
conclude $\V_4(n)$ is also decomposable and the decoding
complexity for one codeword will be $O(N)+O(L^{3/7})$. In
Figure~\ref{figure-three} we compare the performance of $\V_4(2)$
of $103$ elements with a $2$ dimensional constellation of $3$
elements:
$$
\{I_2, A, B\},
$$
where $A=\diag(e^{i2\pi/3}, e^{i2\pi/3})$ and $B=A^2$. It is well
known that this $2$ dimensional constellation has the optimal
diversity product over all the $2$ dimensional unitary
constellations with $3$ elements. Since a large number of
transmitting antennas guarantee that the full diversity at the
transmitter side can be utilized more efficiently (see
Inequality~\ref{pair-err}), it is not too surprising to see that
the $8$ dimensional GPSK constellation performs better at high SNR
region. The figure shows that from $6$db, $8$ dimensional
constellation begins to outperform the $2$ dimensional one.
Around Block Error Rate (BLER) of $10^{-6}$, the performance gain
is about $5$db, which is very pronounced considering the decoding
of $\V_4(2)$ is relatively easy.

\Section{Conclusions and Future Work}

The complex and real orthogonal coding schemes admit simple
decoding algorithms. Based on these schemes, we generalize one
dimensional PSK signals and explicitly construct GPSK unitary
space time constellations. These constellations can be viewed as
higher dimensional generalizations of one dimensional PSK signals
and theoretical analysis shows that their decoding procedures are
decomposable, i.e., the demodulation of these codes can be boiled
down to one dimensional PSK demodulation. Therefore our
constellations have very simple decoding procedures. For some of
the resulting codes (for example, $\V_1(n), \V_2(n)$), the
complexity of ML decoding does not even depend on the transmission
rate. We reallocate the power among the antennas (meanwhile
keeping the total energy) to optimize the diversity product.
Numerical experiments show that our codes perform better than
some of the currently existing comparable ones. For the sample
programs regarding how to construct the proposed constellations,
we refer to~\cite{ha03u2}.

The theoretical analysis~\cite{ta99} shows that the complex and
real orthogonal design constellations only exist on $2, 4$ or $8$
dimension. Based on these schemes we only propose $2, 4$ and $8$
dimensional GPSK constellations. We are trying to generalize the
idea of decomposable GPSK constellations for any dimension.
Constellations with simple decoding algorithms based on the
generalized orthogonal designs~\cite{ja01} are under
investigations.

\vspace{0.3cm}

\textbf{Acknowledgement}

The author thanks J. Rosenthal and V. Tarokh for the helpful
discussions and their insightful suggestions.


\begin{thebibliography}{10}

\bibitem{al98} S.~M. Alamouti.  \newblock A simple transmitter
  diversity scheme for wireless communications.  \newblock {\em
    IEEE J. Selected Areas of Commun.}, pages 1451--1458, October
  1998.

\bibitem{fo96a1} G.~J. Foschini.  \newblock Layered space-time
  architecture for wireless communication in a fading environment
  when using multi-element antennas.  \newblock {\em Bell Labs
    Tech. J.}, 1(2):41--59, 1996.

\bibitem{ha03u2} G.~Han and J.~Rosenthal.  \newblock A website of
  unitary space time constellations with large diversity.
  \newblock \href{http://www.nd.edu/\~{
      }eecoding/space-time/}{http://www.nd.edu/\~{
      }eecoding/space-time/}.

\bibitem{ha02a} B.~Hassibi and B.~M. Hochwald.  \newblock Cayley
  differential unitary space-time codes.  \newblock {\em IEEE
    Trans. Inform. Theory}, 48(6):1485--1503, 2002.  \newblock
  Special issue on Shannon theory: perspective, trends, and
  applications.

\bibitem{ho00} B.~Hochwald and W.~Sweldens.  \newblock
  Differential unitary space-time modulation.  \newblock {\em
    IEEE Trans. Comm.}, pages 2041--2052, December 2000.

\bibitem{ho00a} B.~M. Hochwald and T.~L. Marzetta.  \newblock
  Unitary space-time modulation for multiple-antenna
  communications in {R}ayleigh flat fading.  \newblock {\em IEEE
    Trans. Inform. Theory}, 46(2):543--564, 2000.

\bibitem{ja01} H.~Jafarkhani and V.~Tarokh.  \newblock Multiple
  transmit antenna differential detection from generalized
  orthogonal designs.  \newblock {\em IEEE Trans. Inform.
    Theory}, 47(6):2626--2631, 2001.

\bibitem{ta99} V.~Tarokh, H.~Jafarkhani, and A.~R. Calderbank.
  \newblock Space-time block codes from orthogonal designs.
  \newblock {\em IEEE Trans. Inform. Theory}, 45(5):1456--1467,
  1999.

\bibitem{ta98} V.~Tarokh, N.~Seshadri, and A.~R. Calderbank.
  \newblock Space-time codes for high data rate wireless
  communication: performance criterion and code construction.
  \newblock {\em IEEE Trans. Inform. Theory}, 44(2):744--765,
  1998.

\bibitem{te99} {\.I}.~E. Telatar.  \newblock Capacity of
  multi-antenna {G}aussian channels.  \newblock {\em European
    Trans. Telecommun.}, pages 585--595, 1999.

\end{thebibliography}

\end{document}